\documentclass{article}

\usepackage{amssymb}

\newcommand{\proof}{{\bf Proof:  }}
\newcommand{\remark}{{\bf Remark:  }}
\newcommand{\remarks}{{\bf Remarks:  }}
\newcommand{\example}{{\bf Example:  }}
\newcommand{\examples}{{\bf Examples:  }}

\newcommand{\hb}{\newline\hspace*{\fill}$\Box$}

\newtheorem{theorem}{Theorem}[section]
\newtheorem{lemma}[theorem]{Lemma}
\newtheorem{definition}[theorem]{Definition}
\newtheorem{proposition}[theorem]{Proposition}
\newtheorem{corollary}[theorem]{Corollary}

\begin{document}

\parindent0pt

\title{\bf Cohomology of non-commutative Hilbert schemes}

\author{Markus Reineke\\ Fachbereich Mathematik\\ Bergische Universit\"at Wuppertal\\ Gau\ss str. 20\\ D - 42097 Wuppertal, Germany\\ email: reineke@math.uni-wuppertal.de}

\date{}

\maketitle

\begin{abstract} Non-commutative Hilbert schemes, introduced by M.~V.~Nori, pa\-ra\-met\-rize left ideals of finite codimension in free algebras. More generally, parameter spaces of finite codimensional submodules of free modules over free algebras are considered. Cell decompositions of these varieties are constructed, whose cells are parametrized by certain types of forests. Asymptotics for the corresponding Poincar\'e polynomials and properties of their generating functions are discussed.
\end{abstract}

\section{Introduction}\label{introduction}

\subsection{Non-commutative Hilbert schemes}

Let $k$ be an arbitrary algebraically closed field.
Denote by $A$ the free associative algebra $A=k\langle x_1,\ldots, x_m\rangle$ in $m$ letters.
A representation of $A$ on a finite dimensional $k$-vector space $W$ of dimension $d$ consists of a tuple $\varphi_*$ of $m$ linear endomorphisms $\varphi_1,\ldots,\varphi_m$ of $W$.\\[1ex]
The group ${\rm GL}(W)$ acts on the space ${\rm End}(W)^m$ via base change in $W$, thus $$g(\varphi_1,\ldots,\varphi_m)=(g\varphi_1g^{-1},\ldots,g\varphi_mg^{-1}).$$
The orbits of ${\rm GL}(W)$ in ${\rm End}(W)^m$ correspond bijectively to the isomorphism classes of $d$-dimensional representations of $A$.\\[1ex]
Denote the quotient variety ${\rm End}(W)^m//{\rm GL}(W)$ by $V_d^{(m)}$. Thus, $V_d^{(m)}$ is the affine $k$-variety with coordinate ring isomorphic to the invariant ring $$R=k[{\rm End}(W)^m]^{{\rm GL}(W)}.$$
The $k$-points of $V_d^{(m)}$ are in bijection with the $d$-dimensional {\it semisimple} representations of $A$ (see \cite{Artin}), since an orbit of ${\rm GL}(W)$ in ${\rm End}(W)^m$ is closed if and only if the corresponding representation of $A$ is semisimple.\\[1ex]
By a result of C.~Procesi \cite{Procesi}, the ring $R$ is generated (at least if $k$ is of characteristic $0$) by the functions
$$(\varphi_1,\ldots,\varphi_m)\mapsto {\rm tr}(\varphi_{i_1}\cdots\varphi_{i_s})$$
for sequences $(i_1\ldots i_s)$ in $\{1,\ldots,m\}$. In fact, $R$ is already generated by such functions for $s\leq d^2+1$.\\[1ex]
One of the problems in dealing with the varieties $V_d^{(m)}$ is that they are highly singular, except in the cases $m=1$ (where $V_d^{(m)}\simeq{\bf A}^d$), $d=1$ (where $V_d^{(m)}\simeq{\bf A}^m$) and $d=2=m$ (where $V_d^{(m)}\simeq{\bf A}^5$). Although no explicit desingularizations of the varieties $V_n^{(m)}$ are known (except in case $d=2$; see \cite{Nori}), there exists nevertheless a class of closely related smooth varieties, namely the non-commutative Hilbert schemes introduced in \cite{Nori}. Their construction (in a slightly more general setup) will be reviewed in the following.\\[2ex]
Fix another $k$-vector space $V$ of dimension $n$, together with a basis $v_1,\ldots,v_n$.
Analogous to the above, the affine space $X={\rm Hom}(V,W)\oplus{\rm End}(W)^m$ pa\-ra\-metrizes $d$-dimensional representations of $A$, together with a fixed linear map from $V$ to $W$.
Again, the group ${\rm GL}(W)$ acts on $X$ via base change in $W$, that is, $$g(f,\varphi_1,\ldots,\varphi_m)=(gf,g\varphi_1g^{-1},\ldots,g\varphi_mg^{-1})$$ for $g\in{\rm GL}(W)$, $f\in{\rm Hom}(V,W)$ and $\varphi_1,\ldots,\varphi_m\in{\rm End}(W)$.\\[1ex]
Define a tuple $(f,\varphi_*)$ in $X$ to be stable if $$k\langle \varphi_1,\ldots,\varphi_m\rangle f(V)=W;$$
that is, the image of $f$ generates $W$ as a representation of $A$. Denote by $X^s$ the subset of $X$ consisting of stable tuples. It is easy to see that the stabilizer in ${\rm GL}(W)$ of any stable tuple is trivial.
\begin{definition}[\cite{Nori}] Define ${\rm H}_{d,n}^{(m)}=X^s/{\rm GL}(W)$ as the quotient variety of $X^s$ by ${\rm GL}(W)$.
\end{definition}
The variety ${\rm H}_{d,n}^{(m)}$ is smooth and irreducible, of dimension $N=nd+(m-1)d^2$. It is a principal ${\rm GL}(W)$-bundle over ${\rm H}_{d,n}^{(m)}$.\\[1ex]
All these facts can be easily seen by viewing $X^s$ as the set of stable representations (in the sense of \cite{Ki}) of the quiver $Q_n^{(m)}$ with set of vertices $\{i,j\}$, $n$ arrows from $i$ to $j$, and $m$ loops at $j$, for the dimension vector $e_i+de_j$, and stability $-de_i^*+e_j^*$ (see \cite{Rehns}).\\[1ex]
The variety ${\rm H}_{d,n}^{(m)}$ has several different interpretations:

\begin{lemma}\label{interp} The $k$-points of ${\rm H}_{d,n}^{(m)}$ parametrize each of the following sets:
\begin{enumerate}
\item Equivalence classes of $d$-dimensional representations $W$ of $A$, together with an $n$-tuple
of vectors generating $W$ as a representation of $A$.

\item Equivalence classes of $d$-dimensional representations $W$ of $A$, together with a surjective $A$-homomorphism from the free representation $A^n$ to $W$.

\item $A$-subrepresentations of codimension $d$ of the free representation $A^n$.

\item Isomorphism classes of stable representations of the quiver $Q_n^{(m)}$ as above.
\end{enumerate}
\end{lemma}

\proof The first interpretation follows directly from the definition of ${\rm H}_{d,n}^{(m)}$ as the quotient
of $X^s$ by ${\rm GL}(W)$; the $n$-tuple generating $W$ is $(f(v_1),\ldots,f(v_n))$. The second interpretation
follows by associating to the linear map $f:V\rightarrow W$ the $A$-homomorphism $p_f:A\otimes_kV\rightarrow W$ uniquely
determined by $1\otimes v\mapsto f(v)$ for $v\in V$. Now the tuple $(f,\varphi_*)$ is stable if and only if $p_f$ is
surjective. The third interpretation can be derived from the second: the kernel of $p_f$ is obviously a subrepresentation
of $A\otimes V$ of codimension $d$, and any such subrepresentation $U\subset A\otimes V$ yields a surjection onto a
$d$-dimensional $A$-representation by considering the canonical projection $A\otimes V\rightarrow (A\otimes V)/U$.
The fourth interpretation again follows from the definitions.\hb

\remarks~
\begin{itemize}
\item In particular, the variety ${\rm H}_{d,1}^{(m)}$ parametrizes left ideals of codimension $d$ in $A$. Thus, it can be viewed as a non-commutative Hilbert scheme for the free algebra in $m$ generators, in the same way as the Hilbert scheme ${\rm Hilb}^d({\bf A}^m)$ parametrizes ideals of codimension $d$ in the polynomial ring $k[x_1,\ldots,x_m]$. Under this name, this variety was introduced in \cite{Nori}.
\item Another interpretation of ${\rm H}_{d,1}^{(m)}$ is given in \cite{VDB}, in the framework of Brauer-Severi varieties.
\item The fourth interpretation in Lemma \ref{interp} shows that ${\rm H}_{d,n}^{(m)}$ can also be viewed as a {\it framed} moduli space for representations of the $m$-loop quiver, in the sense of \cite{Nakajima}. See the forthcoming paper \cite{ReFramed} for a general discussion of framed quiver moduli.
\item In the special case $m=1$, the varieties ${\rm H}_{d,n}^{(1)}$ can be described much more explicitly as vector bundles over Grassmannians; see \cite{LR}.
\end{itemize}

Coming back to the original problem of finding a variety closely related to $V_d^{(m)}$, but with better geometric properties, we have the following:\\[1ex]
The obvious map $X^s\rightarrow{\rm End}(W)^m$, forgetting the extra datum $f\in{\rm Hom}(V,W)$, is ${\rm GL}(W)$-equivariant. Thus it descends to a projective morphism $$p:{\rm H}_{d,n}^{(m)}\rightarrow V_d^{(m)}$$ on the level of quotients by ${\rm GL}(W)$. Although the fibres of $p$ are difficult to determine in general, they are at least tractable using the Luna stratification of $V_d^{(m)}$ (see \cite{LBP}) and the theory of nullcones of quiver representations \cite{LB}. This is worked out in \cite{SLB} and allows to derive, for example, a determination of the irreducible components of the fibres and their dimensions.\\[1ex]
It should also be noted that the morphism $p$ extends the canonical map from the Hilbert scheme ${\rm Hilb}^d({\bf A}^m)$ to the $d$-th symmetric power $({\bf A}^m)^d/S_d$. More precisely, we have the following commutative diagram:
$$\begin{array}{ccc}{\rm Hilb}^d({\bf A}^m)&\rightarrow&{\rm H}_{d,1}^{(m)}\\
\downarrow&&\downarrow\\
({\bf A}^m)^d/S_d&\rightarrow&V_d^{(m)}.\end{array}$$
The map in the top row is given by pulling back an ideal via the canonical map $k\langle x_1,\ldots,x_m\rangle \rightarrow k[x_1,\ldots,x_m]$. The map in the bottom row is given by viewing points of ${\bf A}^m$ as one-dimensional representations of $A$, and mapping a $d$-tuple of such to their direct sum.

\subsection{Results}

The main aim of this note is to prove the following results on the geometry of the varieties ${\rm H}_{d,n}^{(m)}$:

\begin{theorem}\label{main1} The variety ${\rm H}_{d,n}^{(m)}$ has a cell decomposition, whose cells are pa\-ra\-me\-trized by $m$-ary forests with $n$ roots and $d$ nodes.
\end{theorem}

The precise definition of such forests (and related combinatorial objects) will be given in section \ref{wordsandforests}. The cell decomposition will be constructed in section \ref{celldec}.\\[1ex]
As a consequence, the Betti numbers in cohomology of ${\rm H}_{d,n}^{(m)}({\bf C})$ can be described. This is possible in a compact form by assembling the Poincar\'e polynomials into a generating function.

\begin{theorem}\label{main2} Define
$$\overline{\zeta}_n^{(m)}(q,t)=\sum_{d=0}^\infty q^{(m-1)\frac{d(d-1)}{2}+(n-1)d}\sum_k\dim H^k({\rm H}_{d,n}^{(m)})q^{-k/2}t^d\in{\bf Q}[q][[t]].$$
Then, as an element of ${\bf Q}[q][[t]]$, the generating function $\overline{\zeta}_n^{(m)}(q,t)$ is determined by the functional equations
$$\overline{\zeta}_n^{(m)}(q,t)=\prod_{i=0}^{n-1}\overline{\zeta}_1^{(m)}(q,q^it)\mbox{ and }\overline{\zeta}_1^{(m)}=1+t\cdot\prod_{i=0}^{m-1}\overline{\zeta}_1^{(m)}(q,q^it).$$
In particular, the cohomological Euler characteristic of ${\rm H}_{d,n}^{(m)}$ is given by
$$\chi({\rm H}_{d,n}^{(m)})=\frac{n}{(m-1)d+n}{{md+n-1}\choose{d}}.$$

\end{theorem}

This theorem will be proved in section \ref{properties}, after computing the Poincar\'e polynomials of the varieties ${\rm H}_{d,n}^{(m)}$ in section \ref{applications}.\\[1ex]
Finally, in section \ref{asymp}, we consider the asymptotic behaviour of both the Euler characteristic and the Poincar\'e polynomials, in the spirit of \cite{Kon}. The main result is:

\begin{theorem} After a suitable renormalization, the distribution of the Betti numbers of ${\rm H}_{d,1}^{(m)}$ for large $d$ has the Airy distribution as a limit law.
\end{theorem}

The results of Theorem \ref{main2} were predicted by computer experiments using the formula \cite[Theorem 6.7]{Rehns} for the cohomology of quiver moduli. Namely, these experiments suggested the above formula for the Euler characteristic.
Since these numbers are just the Catalan numbers in the special case $m=2$, $n=1$, a relation to trees could be expected.
This relation could then be established by using a construction in \cite{VDB}; see section \ref{celldec}.\\[2ex]
{\bf Acknowledgment:} This paper was written during a stay at the University of Antwerp, with the aid of a grant of the European Science Foundation in the frame of the Priority Programme ``Noncommutative Geometry". I would like to thank L.~Le Bruyn for pointing out to me the role of non-commutative Hilbert schemes, in particular the works \cite{Nori,VDB}, and for many interesting discussions. I am grateful to P.~Biane for pointing me to the probabilistic literature, and to
P.~Duchon for advice on how to use his results.

\section{Words and forests}\label{wordsandforests}

Denote by $\Omega=\Omega^{(m)}$ the set of finite words in the alphabet $\{1,\ldots,m\}$. The length $l(w)$ of a word $w=(i_1\ldots i_s)$ is defined to be $s$. The unique word of length $0$ is denoted by $()$. The concatenation of words $w$ and $w'$ is written as $ww'$; in particular, the concatenation of a word $w$ with a letter $i$ is denoted by $wi$.\\[1ex]
For a word $w=(i_1\ldots i_s)\in\Omega$, define $x_\omega=x_{i_s}\ldots x_{i_1}\in A$. The elements $x_w$ for $w\in\Omega$ thus form a $k$-basis for $A$. Furthermore, for any representation of $A$ given by a tuple $(\varphi_1,\ldots,\varphi_m)$ of endomorphisms of the vector space $W$, define $\varphi_w=\varphi_{i_s}\cdots\varphi_{i_1}\in{\rm End}_k(W)$.\\[1ex]
The set $\Omega$ can also be viewed as the free $m$-ary tree, with the empty word $()$ as the root, and with each word $w$ having $w1,\ldots,wm$ as its $m$ successors.\\[1ex]
\example The set $\Omega^{(3)}$ can be viewed as the free ternary tree as follows:

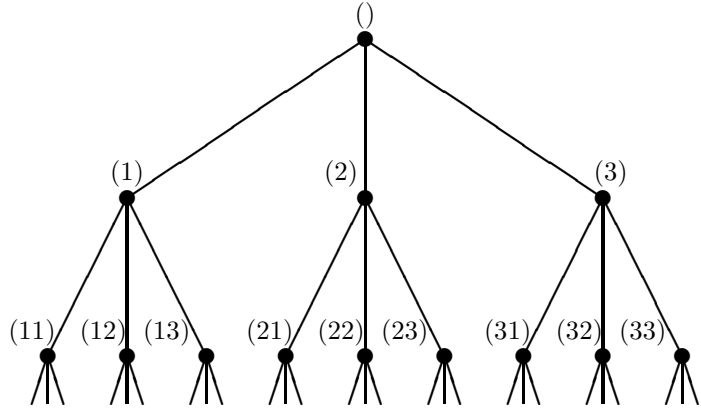
\begin{center}\label{fig:a1}
\setlength{\unitlength}{3pt}
\begin{picture}(80,60)(0,-10)
\thicklines
\multiput(0,0)(30,0){3}{\line(1,2){10}}
\multiput(10,0)(30,0){3}{\line(0,1){20}}
\multiput(20,0)(30,0){3}{\line(-1,2){10}}
\multiput(0,0)(10,0){9}{\circle*{2}}
\multiput(10,20)(30,0){3}{\circle*{2}}
\put(10,20){\line(3,2){30}}
\put(40,20){\line(0,1){20}}
\put(70,20){\line(-3,2){30}}
\put(40,40){\circle*{2}}
\put(40,43){\makebox(0,0){$()$}}
\put(10,23){\makebox(0,0){$(1)$}}
\put(37,23){\makebox(0,0){$(2)$}}
\put(71,23){\makebox(0,0){$(3)$}}
\put(-2,3){\makebox(0,0){$(11)$}}
\put(7,3){\makebox(0,0){$(12)$}}
\put(15,3){\makebox(0,0){$(13)$}}
\put(28,3){\makebox(0,0){$(21)$}}
\put(37,3){\makebox(0,0){$(22)$}}
\put(45,3){\makebox(0,0){$(23)$}}
\put(58,3){\makebox(0,0){$(31)$}}
\put(67,3){\makebox(0,0){$(32)$}}
\put(75,3){\makebox(0,0){$(33)$}}
\multiput(0,0)(10,0){9}{\line(-1,-3){2}}
\multiput(0,0)(10,0){9}{\line(0,-1){6}}
\multiput(0,0)(10,0){9}{\line(1,-3){2}}

\end{picture}
\end{center}

Consequently, we view $\{1,\ldots,n\}\times\Omega$ as the free $m$-ary forest with $n$ roots, since by definition a forest is nothing else than an ordered tuple of trees.\\[1ex]
Since any finite $m$-ary tree can be embedded into the free tree $\Omega$, we arrive at the following definition:
\begin{definition}~
\begin{enumerate}
\item A finite subset $S\subset\Omega$ is called an $m$-ary tree if it is closed under taking left subwords, that is, $w\in S$ provided $ww'\in S$ for some $w'\in\Omega$. 
\item An $m$-ary forest with $n$ roots is an $n$-tuple of (possibly empty) $m$-ary trees $S_*=(S_1,\ldots,S_n)$.
\end{enumerate}
\end{definition}

Via this definition, the nodes of a forest $S_*$ are parametrized by pairs $(k,w)$, where $w\in S_k$ and $k=1,\ldots,n$. We simply write $(k,w)\in S_*$. In particular, the number of nodes $|S_*|$ of a forest $S_*$ is just the sum of the cardinalities
$$|S_*|=\sum_{i=1}^n|S_i|.$$
Denote the set of $m$-ary forests with $n$ roots and $d$ nodes by ${\cal F}_{d,n}^{(m)}$.\\[1ex]
A pair $(k,w)\in\{1,\ldots,m\}\times\Omega$ is called critical for a forest $S_*$ if either $w=()$ and $S_k=\emptyset$,
or $w\not\in S_k$, but $w'\in S_k$ for $w=w'i$. Denote by $C(S_*)$ the set of pairs $(k,w)$ which are critical for $S_*$,
and denote the cardinality of the critical set $C(S_*)$ by $c(S)$. An easy induction on the cardinality of the forest
$S_*$ immediately yields the formula
$$c(S_*)=(m-1)|S_*|+n.$$
The set $\Omega$ is totally ordered with respect to the lexicographical ordering:
suppose $w=(i_1\ldots i_s)$ and $w'=(j_1\ldots j_t)$. Let $k$ be the maximal index such that $i_k=j_k$
(we set $k=0$ if no such index exists). Then $w\leq w'$ if either $k=s=l(w)$, or $i_{k+1}\leq j_{k+1}$.\\[1ex]
We extend this to a total ordering on trees by setting $S\leq S'$ if either $|S|>|S'|$, or $|S|=|S'|$ and, writing
$S=\{w_1<\ldots<w_s\}$ and $S'=\{w_1'<\ldots<w_s'\}$, we have $w_p<w_p'$ for the minimal $p$ such that $w_p\not=w_p'$.\\[1ex]
Finally, we extend this to a total ordering on each set ${\cal F}_{d,n}^{(m)}$ by defining $S_*\leq S_*'$ if $S_*=S_*'$ or
$S_p<S_p'$ for the minimal $p$ such that $S_p\not=S_p'$.\\[1ex]
We also order the free forest $\{1,\ldots,n\}\times\Omega$ by setting $(k,w)\leq(l,w')$ if either $k<l$, or $k=l$ and
$w\leq w'$. In particular, this induces a total ordering on the nodes of any forest $S_*$.\\[2ex]
\remark If $S<S'$ are trees of the same cardinality, then $w_p\not\in S'$ for the element $w_p$ defined above: otherwise,
$w_p=w_q'$ for some $q<p$ since $w_p<w_p'$. But then $w_p=w_q'=w_q$ by the choice of $p$; a contradiction.\\[2ex]
For each forest $S_*$, define $D(S_*)$ as the set of all tuples $(k,w,l,w')$ in
$$(\{1,\ldots,m\}\times\Omega)\times(\{1,\ldots,m\}\times\Omega)$$ such that
$$(k,w)\in S_*,\;\; (l,w')\in C(S_*),\;\; (k,w)<(l,w').$$
Denote by $d(S_*)$ the cardinality $|D(S_*)|$ of $D(S_*)$.\\[1ex]
\example The binary trees with four nodes are given in their lexicographical ordering as follows:

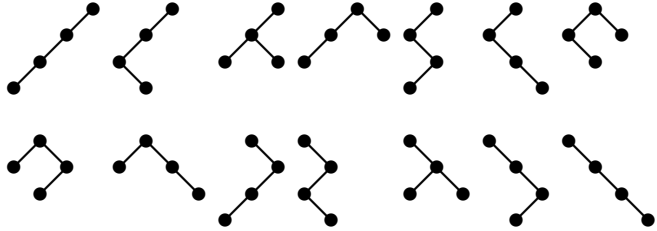
\begin{center}
\setlength{\unitlength}{10pt}
\begin{picture}(24,8)(0,0)
\thicklines
\put(0,5){\circle*{.5}}
\put(1,6){\circle*{.5}}
\put(2,7){\circle*{.5}}
\put(3,8){\circle*{.5}}

\put(4,6){\circle*{.5}}
\put(5,5){\circle*{.5}}
\put(5,7){\circle*{.5}}
\put(6,8){\circle*{.5}}

\put(8,6){\circle*{.5}}
\put(9,7){\circle*{.5}}
\put(10,6){\circle*{.5}}
\put(10,8){\circle*{.5}}

\put(11,6){\circle*{.5}}
\put(12,7){\circle*{.5}}
\put(13,8){\circle*{.5}}
\put(14,7){\circle*{.5}}

\put(15,5){\circle*{.5}}
\put(16,6){\circle*{.5}}
\put(15,7){\circle*{.5}}
\put(16,8){\circle*{.5}}

\put(18,7){\circle*{.5}}
\put(19,8){\circle*{.5}}
\put(19,6){\circle*{.5}}
\put(20,5){\circle*{.5}}

\put(21,7){\circle*{.5}}
\put(22,8){\circle*{.5}}
\put(22,6){\circle*{.5}}
\put(23,7){\circle*{.5}}

\put(0,2){\circle*{.5}}
\put(1,3){\circle*{.5}}
\put(1,1){\circle*{.5}}
\put(2,2){\circle*{.5}}

\put(4,2){\circle*{.5}}
\put(5,3){\circle*{.5}}
\put(6,2){\circle*{.5}}
\put(7,1){\circle*{.5}}

\put(8,0){\circle*{.5}}
\put(9,1){\circle*{.5}}
\put(10,2){\circle*{.5}}
\put(9,3){\circle*{.5}}

\put(11,3){\circle*{.5}}
\put(12,2){\circle*{.5}}
\put(11,1){\circle*{.5}}
\put(12,0){\circle*{.5}}

\put(15,3){\circle*{.5}}
\put(15,1){\circle*{.5}}
\put(16,2){\circle*{.5}}
\put(17,1){\circle*{.5}}

\put(18,3){\circle*{.5}}
\put(19,2){\circle*{.5}}
\put(20,1){\circle*{.5}}
\put(19,0){\circle*{.5}}

\put(21,3){\circle*{.5}}
\put(22,2){\circle*{.5}}
\put(23,1){\circle*{.5}}
\put(24,0){\circle*{.5}}

\put(0,5){\line(1,1){1}}
\put(1,6){\line(1,1){1}}
\put(2,7){\line(1,1){1}}
\put(4,6){\line(1,1){1}}
\put(5,7){\line(1,1){1}}
\put(8,6){\line(1,1){1}}
\put(9,7){\line(1,1){1}}
\put(11,6){\line(1,1){1}}
\put(12,7){\line(1,1){1}}
\put(15,7){\line(1,1){1}}
\put(15,5){\line(1,1){1}}
\put(18,7){\line(1,1){1}}
\put(21,7){\line(1,1){1}}
\put(0,2){\line(1,1){1}}
\put(4,2){\line(1,1){1}}
\put(1,1){\line(1,1){1}}
\put(8,0){\line(1,1){1}}
\put(9,1){\line(1,1){1}}
\put(11,1){\line(1,1){1}}
\put(15,1){\line(1,1){1}}
\put(19,0){\line(1,1){1}}

\put(4,6){\line(1,-1){1}}
\put(9,7){\line(1,-1){1}}
\put(13,8){\line(1,-1){1}}
\put(15,7){\line(1,-1){1}}
\put(18,7){\line(1,-1){1}}
\put(19,6){\line(1,-1){1}}
\put(21,7){\line(1,-1){1}}
\put(22,8){\line(1,-1){1}}
\put(1,3){\line(1,-1){1}}
\put(5,3){\line(1,-1){1}}
\put(6,2){\line(1,-1){1}}
\put(9,3){\line(1,-1){1}}
\put(11,3){\line(1,-1){1}}
\put(11,1){\line(1,-1){1}}
\put(15,3){\line(1,-1){1}}
\put(16,2){\line(1,-1){1}}
\put(18,3){\line(1,-1){1}}
\put(19,2){\line(1,-1){1}}
\put(21,3){\line(1,-1){1}}
\put(22,2){\line(1,-1){1}}
\put(23,1){\line(1,-1){1}}

\end{picture}
\end{center}

\section{The cell decomposition}\label{celldec}

The starting point for the construction of a cell decomposition of ${\rm H}_{d,n}^{(m)}$ is the construction, due to M.~Van den Bergh \cite{VDB}, of a covering by open affine set for the varieties ${\rm H}_{d,1}^{(m)}$. We will first translate this to the language of forests. Then, by using the total ordering on forests, successive complements of these open affine sets will be shown to yield a cell decomposition. This might be compared with the construction of the cell decomposition by Schubert cells for Grassmannians (which constitute the case $m=0$ in the present setup).\\[1ex]
Note that, in contrast to the commutative Hilbert schemes ${\rm Hilb}^d({\bf A}^2)$ (see \cite{ES}), it is not possible to apply the Bialynicki-Birula method \cite{BiaBi} to the varieties ${\rm H}_{d,n}^{(m)}$, since the natural torus action has infinitely many fixed points; see the example at the end of this section.

\begin{definition} For each forest $S_*$ in ${\cal F}_{d,n}^{(m)}$, define $U_{S_*}$ as the set of all tuples $(f,\varphi_*)\in {\rm H}_{d,n}^{(m)}$ such that the set
$$\{\varphi_wf(v_k)\, :\, (k,w)\in S_*\}$$
forms a basis of $W$.
\end{definition}

\remark From the definition of $U_{S_*}$, it is clear that $\dim k\langle \varphi_*\rangle f(v_k)\geq|S_k|$ for all $k=1,\ldots,m$, and all $(f,\varphi_*)\in U_{S_*}$.\\[2ex]
Let $(f,\varphi_*)$ be a point in $U_{S_*}$. For any $(l,w')\in C(S_*)$, we can write uniquely $$\varphi_{w'}f(v_l)=\sum_{(k,w)\in S_*}\lambda_{(k,w,l,w')}\varphi_wf(v_k).$$
It is easy to see that the map assigning to $(f,\varphi_*)$ the coefficient $\lambda_{(k,w,l,w')}$ induces an algebraic function $$\Lambda_{(k,w,l,w')}:U_{S_*}\rightarrow k.$$

\begin{lemma}\label{bel} Let $(f,\varphi_*)$ be a point in ${\rm H}_{d,n}^{(m)}$. Let $\overline{S}_*$ be a forest such that the vectors $\varphi_wf(v_k)$ for $(k,w)\in\overline{S}_*$ are linearly independent in $W$. Then there exists a forest $S_*\in{\cal F}_{d,n}^{(m)}$, containing $\overline{S}_*$, such that $(f,\varphi_*)\in U_{S_*}$.
\end{lemma}

\proof We proceed by downward induction on the cardinality $|\overline{S}_*|$ of $\overline{S}_*$.
If this equals $d$, there is nothing to prove.
Otherwise, let $U$ be the span of the vectors
$\varphi_wf(v_k)$ for $(k,w)\in\overline{S}_*$,
and consider $k\langle \varphi_*\rangle U\subset W$.\\[1ex]
If this space strictly contains $U$, we can choose an
index $k$, a word $w\in\overline{S}_k$ and some $i\in\{1,\ldots,m\}$ such that
$\varphi_{wi}f(v_k)\not\in U$. Define $\overline{S}_l'$ as $\overline{S}_l$ if $l\not=k$,
and as $\overline{S}_k\cup\{wi\}$ if $l=k$. By definition, $\overline{S}_*'$ forms a forest
fulfilling the conditions of the lemma, of strictly larger cardinality than $\overline{S}_*$.
By induction, we are done.\\[1ex]
If $k\langle\varphi_*\rangle U=U$, we have $f(V)\not\subset U$, since
otherwise $W=k\langle\varphi_*\rangle f(V)\subset k\langle\varphi_*\rangle U=U\not=W$, a contradiction.
Thus, there exists an index $k$ such that $\overline{S}_k=\emptyset$ and $f(v_k)\not\in U$.
Defining $\overline{S}_l'$ as $\overline{S}_l$ for $l\not=k$, and as $\{()\}$ for $l=k$,
we can proceed as above. \hb

Applying the lemma to the empty forest, we get:

\begin{corollary} The variety ${\rm H}_{d,n}^{(m)}$ is the union of the $U_{S_*}$, ranging over all forests $S_*$ in ${\cal F}_{d,n}^{(m)}$.
\end{corollary}

\begin{lemma} Each $U_{S_*}$ is an open subset of ${\rm H}_{d,n}^{(m)}$, isomorphic to an affine space of dimension $N=(m-1)d^2+nd$.
\end{lemma}

\proof The defining condition of $U_{S_*}$ can be rephrased as the non-vanishing of a determinant, and is therefore open. The isomorphism to an affine space is easily seen to be provided by the functions $\Lambda_{(k,w,l,w')}$ on $U_{S_*}$ for $(k,w)\in S_*$ and $(l,w')\in C(S_*)$: note that, by the formula for $c(S_*)$ above, there are precisely $|S_*|\cdot((m-1)|S_*|+n)$ of them.\hb

Now we come to the definition of certain locally closed subvarieties of ${\rm H}_{d,n}^{(m)}$, which will be shown to provide a cell decomposition.

\begin{definition} For each forest $S_*$, define $Z_{S_*}$ as the set of all points $(f,\varphi_*)$ of $U_{S_*}$ such that the following holds:\\[1ex]
For all critical pairs $(k,w)\in C(S_*)$, the vector $\varphi_wf(v_k)$ is contained in the span of the vectors $\varphi_{w'}f(v_l)$ for $(l,w')<(k,w)$.
\end{definition}

\example Consider the case $m=3$, $d=6$, $n=3$. Let $S_*$ be the following forest:

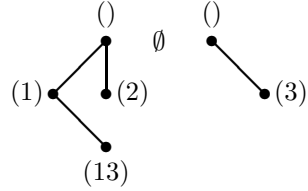
\begin{center}
\setlength{\unitlength}{2pt}
\begin{picture}(40,35)(0,-10)
\thicklines
\put(0,10){\circle*{2}}
\put(-5,10){\makebox(0,0){$(1)$}}
\put(10,0){\circle*{2}}
\put(10,-5){\makebox(0,0){$(13)$}}
\put(10,20){\circle*{2}}
\put(10,25){\makebox(0,0){$()$}}
\put(10,10){\circle*{2}}
\put(15,10){\makebox(0,0){$(2)$}}
\put(20,20){\makebox(0,0){$\emptyset$}}
\put(30,20){\circle*{2}}
\put(30,25){\makebox(0,0){$()$}}
\put(40,10){\circle*{2}}
\put(45,10){\makebox(0,0){$(3)$}}
\put(10,20){\line(-1,-1){10}}
\put(0,10){\line(1,-1){10}}
\put(10,20){\line(0,-1){10}}
\put(30,20){\line(1,-1){10}}

\end{picture}
\end{center}

In the notation introduced above, we have
$$S_*=(\{(),(1),(13),(2)\},\emptyset,\{(),(3)\}).$$
The critical set $C(S_*)$ consists of the pairs
\begin{eqnarray*}C(S_*)&=&\{(1,(11)),(1,(12)),(1,(131)),(1,(132)),(1,(133)),\\
&&(1,(21)),(1,(22)),(1,(23)),(1,(3)),(2,()),\\
&&(3,(1)),(3,(2)),(3,(31)),(3,(32)),(3,(33))\}.\end{eqnarray*}
Assume that $(f,\varphi_*)\in U_{S_*}$. Then the following elements form a basis $B$ of $W$:
$$f(v_1),\,\varphi_1f(v_1),\,\varphi_3\varphi_1f(v_1),\,\varphi_2f(v_1),\,f(v_3),\,\varphi_3f(v_3).$$
If $(f,\varphi_*)\in Z_{S_*}$, then, with respect to the above basis $B$, the matrices representing $f$ and the $\varphi_i$, respectively, are of the form:

$$f=\left[\begin{array}{ccc}1&\Box&0\\ 0&\Box&0\\ 0&\Box&0\\ 0&\Box&0\\ 0&0&1\\ 0&0&0\end{array}\right],\;\;\; \varphi_1=\left[\begin{array}{cccccc}0&\Box&\Box&\Box&\Box&\Box\\ 1&\Box&\Box&\Box&\Box&\Box\\ 0&0&\Box&\Box&\Box&\Box\\ 0&0&0&\Box&\Box&\Box\\ 0&0&0&0&\Box&\Box\\ 0&0&0&0&0&\Box\end{array}\right]$$
$$\varphi_2=\left[\begin{array}{cccccc}0&\Box&\Box&\Box&\Box&\Box\\ 0&\Box&\Box&\Box&\Box&\Box\\ 0&0&\Box&\Box&\Box&\Box\\ 1&0&0&\Box&\Box&\Box\\ 0&0&0&0&\Box&\Box\\ 0&0&0&0&0&\Box\end{array}\right],\;\;\; \varphi_3=\left[\begin{array}{cccccc}\Box&0&\Box&\Box&0&\Box\\ \Box&0&\Box&\Box&0&\Box\\ \Box&1&\Box&\Box&0&\Box\\ \Box&0&0&\Box&0&\Box\\ 0&0&0&0&0&\Box\\ 0&0&0&0&1&\Box\end{array}\right],$$
where $\Box$ indicates an arbitrary entry.\\[2ex]
We now study the properties of the sets $Z_{S_*}$.\\[1ex]
\remark From the definition of $Z_{S_*}$, it is clear that $$\dim \sum_{l=1}^kk\langle \varphi_*\rangle 
f(v_l)=|S_1|+\ldots+|S_k|$$ for all $k=1,\ldots,m$, and all $(f,\varphi_*)\in Z_{S_*}$.

\begin{theorem}\label{cd} For all forests $S_*$, we have $Z_{S_*}=U_{S_*}\setminus\bigcup_{S_*'<S_*}U_{S_*'}$.
\end{theorem}

The proof of this theorem will be given in the two obvious steps; namely, we prove each inclusion in the claimed equality separately.

\begin{lemma} For all forests $S_*$, we have $Z_{S_*}\subset U_{S_*}\setminus\bigcup_{S_*'<S_*}U_{S_*'}$.
\end{lemma}

\proof Let $(f,\varphi_*)$ be a point in $Z_{S_*}$, and suppose that $(f,\varphi_*)$ belongs to
$U_{S_*'}$ for some $S_*'<S_*$. Let $p$ be minimal such that $S_p'\not=S_p$, so that $S_p'<S_p$.
Then, by the remarks following the definitions of $U_{S_*}$ and $Z_{S_*}$, respectively, 
we have $$|S_1|+\ldots+|S_p|=\dim \sum_{k=1}^pk\langle\varphi_*\rangle f(v_k)
\geq|S_1'|+\ldots+|S_p'|.$$
Thus $|S_p'|=|S_p|$ by the choice of $p$.
Writing $$S_p=\{w_1<\ldots< w_s\},\;\;\; S_p'=\{w_1'<\ldots<w_s'\},$$
and choosing $q$ minimal such
that $w_q\not=w_q'$, we thus have $w_q'<w_q$. By the remark following the definition of the
ordering on trees, we get $w_q'\not\in S_p$. Note that this implies $w_q'\not=()$,
since otherwise $()\not\in S_p$, thus $S_p=\emptyset$, contradicting $S_p\not=S_p'$ and
$|S_p|=|S_p'|$. Thus we can write $w_q'=wi$. Since $S_p$ is a tree, we have $w\in S_p'$,
and also $w<w_q'$, thus $w=w_r'=w_r$ for some $r<q$, and thus $w\in S_p$.
\\[1ex]
We arrive at the
situation $w\in S_p$, $wi\not\in S_p$, which by definition of $Z_{S_*}$ means that
$\varphi_{wi}f(v_p)$ is contained in the span of the vectors $\varphi_{w'}f(v_s)$ for $(s,w')<(p,wi)$. If $s<p$, then $w'\in S_s'$ since $S_s'=S_s$. If $s=p$, then $w'<wi=w_q'<w_q$, and thus $w'=w_t=w_t'\in S_p'$ for some $t<q$.\\[1ex]
We conclude that $\varphi_{wi}f(v_p)$ is contained in the span of vectors $\varphi_w'f(v_s)$ for $(s,w')\in S_*'$, which is impossible, since by the assumption $(f,\varphi_*)\in U_{S_*'}$, all the vectors $\varphi_w'f(v_s)$ for $w'\in S_s$ are linearly independent. \hb

\begin{lemma} For all forests $S_*$, we have $Z_{S_*}\supset U_{S_*}\setminus\bigcup_{S_*'<S_*}U_{S_*'}$.
\end{lemma}

\proof Let $(f,\varphi_*)$ be in $U_{S_*}$, and assume that $(f,\varphi_*)\not\in Z_{S_*}$.
Then there exists some pair $(k,w)\in C(S_*)$, such that $\varphi_wf(v_k)$ is not contained in the
span of the vectors $\varphi_{w'}f(v_l)$ for $(l,w')<(k,w)$.\\[1ex]
Among such pairs, choose one with minimal index $k$, and with minimal $w$. Define a forest $\overline{S}_*$ by setting $\overline{S}_l=S_l$ for $l<k$, $\overline{S}_k=\{w'\in S_k\, :\, w'\leq w\}$, and $\overline{S}_l=\emptyset$ for $l>k$. The assumptions of Lemma \ref{bel} are satisfied, yielding a forest $S_*'\supset\overline{S}_*$ such that $(f,\varphi_*)\in U_{S_*'}$. It remains to prove that $S_*'<S_*$, yielding a contradiction.\\[1ex]
We first prove $S_l'=S_l$ for $l<k$: suppose not, then we can choose a word $w'\in S_l'\setminus\overline{S}_l$. Without loss of generality, $w'$ can be chosen to be critical for $\overline{S}_l=S_l$. Since $(f,\varphi_*)\in U_{S_*'}$, the vector $\varphi_{w'}f(v_l)$ is linearly independent of the vectors $\varphi_{w''}f(v_p)$ for $(p,w'')<(l,w')$. But $l<k$, contradicting the minimality of $k$.\\[1ex]
Now we compare $S_k$ and $S_k'$: since $w\not\in S_k$, but $w\in\overline{S}_k\subset S_k'$, we have $S_k'\not=S_k$. By the remark following the definition of $U_{S_*}$, we have $\dim k\langle\varphi_*\rangle f(v_k)\geq|S_k|$.\\[1ex]
We can thus assume that the extension $S_k'$ of $\overline{S}_k$ is chosen in such a way that $|S_k'|\geq|S_k|$. If $|S_k'|>|S_k|$, we have proved $S_*'<S_*$ by definition; so assume $|S_k'|=|S_k|$. We have $w\not=()$, since otherwise $S_k=\emptyset$, and thus $\emptyset=S_k'\supset \overline{S}_k\not=\emptyset$, a contradiction.\\[1ex]
Write $$\overline{S}_k\setminus\{w\}=\{w_1<\ldots<w_{p-1}\}.$$
By definition of $\overline{S}_k$, we then have $$S_k=\{w_1<\ldots<w_{p-1}<w_p<\ldots\}.$$
Assume $w'\in S_k'\setminus\overline{S}_k$. Then $w'>w$, since otherwise, we get a contradiction to the minimality of $w$ as above. Thus, we can write
$$S_k'=\{w_1<\ldots<w_{p-1}<w<\ldots\}.$$
But $w_p\in S_k\setminus\overline{S}_k$, thus $w_p>w$ by definition of $\overline{S}_k$, which proves $S_k'<S_k$. \hb

Thus, Theorem \ref{cd} is proved.

\begin{lemma}\label{aff} The variety $Z_{S_*}$ is isomorphic to an affine space of dimension $d(S_*)$.
\end{lemma}

\proof Identify the affine space ${\bf A}^{d(S_*)}$ with the $k$-space $Y_{S_*}$ with basis $e_{(k,w,l,w')}$ for $(k,w,l,w')\in D(S_*)$. We define a morphism $\Psi_{S_*}:Y_{S_*}\rightarrow Z_{S_*}$. Let a vector $$y=\sum_{(k,w,l,w')\in D(S_*)}\lambda_{(k,w,l,w')}e_{(k,w,l,w')}$$ in $Y_{S_*}$ be given. Choose a basis $\{b_{(k,w)}\}$ of $W$, parametrized by the pairs $(k,w)\in S_*$. Define a point $\Psi_{S_*}(y)$ of $Z_{S_*}$ by the following formulas:
\begin{itemize}
\item $f(v_k)=b_{(k,())}$ if $S_k\not=\emptyset$,
\item $f(v_k)=\sum_{(l,w,k,())\in D(S_*)}\lambda_{(l,w,k,())}b_{(l,w)}$ if $S_k=\emptyset$,
\item $\varphi_i(b_{(k,w)})=b_{(k,wi)}$ if $w\in S_k$, $wi\in S_k$,
\item $\varphi_i(b_{(k,w)})=\sum_{(l,w',k,wi)}\lambda_{(l,w',k,wi)}b_{(l,w')}$ if $w\in S_k$, $wi\not\in S_k$.
\end{itemize}
Using the functions $\Lambda_{(k,w,l,w')}$ for $(k,w,l,w')\in D(S_*)$, the map $\Psi_{S_*}$ thus defined is easily seen to provide an isomorphism $Y_{S_*}\simeq Z_{S_*}$ by definition of $Z_{S_*}$. \hb

{\bf Proof (of Theorem \ref{main1}):} For each forest $S_*$, define $$A_S={\rm H}_{d,n}^{(m)}\setminus\bigcup_{S_*'<S_*}S_*',$$
which is a closed subvariety of ${\rm H}^{(m)}_{d,n}$. Enumerate the set ${\cal F}_{d,n}^{(m)}$ as $S_*^1<\ldots<S_*^u$, and define $A_i=A_{S_*^i}$. Then ${\rm H}_{d,n}^{(m)}=A_1\supset\ldots\supset A_u\supset 0$ defines a filtration of ${\rm H}_{d,n}^{(m)}$ by closed subvarieties. Each of the successive complements $$A_i\setminus A_{i+1}=Z_{S_*^i}$$ is isomorphic to an affine space by Lemma \ref{aff}. \hb

\example  As noted at the beginning of this section, it is not possible to apply the Bialynicki-Birula method to the varieties ${\rm H}_{d,n}^{(m)}$: they carry a natural action of the $m$-torus $(k^*)^m$ via
$$(t_1,\ldots,t_n)(f,\varphi_1,\ldots,\varphi_m)=(f,t_1\varphi_1,\ldots,t_m\varphi_m),$$
but this action has infinitely many fixed points already in the case $n=1$, $m=2$ and $d=4$, as the example
$$f=\left[\begin{array}{c}1\\ 0\\ 0\\ 0\end{array}\right],\;\;\; \varphi_1=\left[\begin{array}{cccc} 0&0&0&0\\ 1&0&0&0\\ 0&0&0&\lambda\\ 0&0&0&0\end{array}\right],\;\;\; \varphi_2=\left[\begin{array}{cccc}0&0&0&0\\ 0&0&0&0\\ 0&1&0&0\\ 1&0&0&0\end{array}\right]$$
with $\lambda\in k$ shows by a direct calculation. Note that $\varphi_1\varphi_2=\lambda\varphi_2\varphi_1$ in this example, which shows that the problem arises precisely due to the non-commutativity of the situation.

\section{Applications of the cell decomposition}\label{applications}

\subsection{Normal forms for representations and submodules}

As an application of the cell decomposition constructed above, we can derive normal forms for representations of $A$ equipped with generating vectors, as well as for finite codimensional subspaces of free modules over free algebras.\\[1ex]
In the first case, we just have to rephrase the proof of Lemma \ref{aff}.
\begin{proposition} If a tuple $(f,\varphi_*)$ belongs to $Z_{S_*}$, then there exists a basis 
$\{b_{(k,w)}\}$ of $W$, parametrized by the pairs $(k,w)\in S_*$, such that $(f,\varphi_*)$ can be written in the form
:\begin{itemize}
\item $f(v_k)=b_{(k,())}$ if $S_k\not=\emptyset$,
\item $f(v_k)=\sum_{(l,w,k,())\in D(S_*)}\lambda_{(l,w,k,())}b_{(l,w)}$ if $S_k=\emptyset$,
\item $\varphi_i(b_{(k,w)})=b_{(k,wi)}$ if $w\in S_k$, $wi\in S_k$,
\item $\varphi_i(b_{(k,w)})=\sum_{(l,w',k,wi)\in D(S_*)}\lambda_{(l,w',k,wi)}b_{(l,w')}$ if $w\in S_k$, $wi\not\in S_k$.
\end{itemize}
\end{proposition}

In the second case, we use the third interpretation of the variety ${\rm H}_{d,n}^{(m)}$ from Lemma \ref{interp}. Recall that, to any tuple $(f,\varphi_*)\in{\rm H}_{d,n}^{(m)}$, we associate the subspace $$U={\rm Ker}(A\otimes V\rightarrow W)\subset A\otimes V.$$
The next result now follows directly from the definition of the cell $Z_{S_*}$.

\begin{proposition} If $U\subset A\otimes V$ corresponds to a tuple
$(f,\varphi_*)\in Z_{S_*}$ for a forest $S_*\in{\cal F}_{d,n}^{(m)}$, then $U$ is generated by elements
$$x_{w'}\otimes v_l-\sum_{(k,w)\, :\, (k,w,l,w')\in D(S_*)}\lambda_{(k,w,l,w')}x_w\otimes v_k,$$
where $(l,w')$ runs over all elements of $C(S_*)$.
\end{proposition}

This gives a normal form for arbitrary finite codimensional subspaces, which - in the case $n=1$ - one might view as a more precise version of \cite[Theorem 6.8]{Cohn}.

\subsection{(Co-)homology}
Applying \cite[Chapter 1]{Fu} to Theorem \ref{main1}, we get immediately:

\begin{corollary}\label{it} The intersection theory $A_*({\rm H}_{d,n}^{(m)})$ is a free abelian group, with a basis given by the classes of the closures $\overline{Z_{S_*}}$, for $S_*$ running through ${\cal F}_{d,n}^{(m)}$.
\end{corollary}

The Poincar\'e polynomial in intersection theory of ${\rm H}_{d,n}^{(m)}$ is thus given as follows:
$$\sum_{i}\dim_{\bf Z}A_i({\rm H}_{d,n}^{(m)})q^i=\sum_{S_*\in{\cal F}_{d,n}^{(m)}}q^{d(S_*)}.$$

{\it From now on, consider the special case $k={\bf C}$ of the complex numbers as the ground field.}\\[1ex]
Then, by \cite[Example 19.1.11.(d)]{Fu}, the cycle map $$A_*({\rm H}_{d,n}^{(m)})\rightarrow H_{2*}^{BM}({\rm H}_{d,n}^{(m)})$$ induces an isomorphism between intersection theory and even Borel-Moore homology (with integer coefficients), whereas the odd Borel-Moore homology vanishes. Moreover, the smoothness of ${\rm H}_{d,n}^{(m)}$ allows us to apply Poincar\'e duality to identify the Borel-Moore homology with the cohomology $$H_k^{BM}({\rm H}_{d,n}^{(m)})\simeq H^{2N-k}({\rm H}_{d,n}^{(m)}),$$
where $N$ equals the dimension of ${\rm H}_{d,n}^{(m)}$. This discussion yields the following formula:

\begin{corollary}\label{betti} The Poincar\'e polynomial in cohomology of ${\rm H}_{d,n}^{(m)}$ is given by

$$\sum_i\dim H^i({\rm H}_{d,n}^{(m)})t^i=\sum_{S_*\in{\cal F}_{d,n}^{(m)}}t^{2((m-1)d^2+nd-d(S_*))}.$$
\end{corollary}

In particular, specializing $q$ to $-1$, we see that the cohomological Euler characteristic of ${\rm H}_{d,n}^{(m)}$ equals the cardinality of ${\cal F}_{d,n}^{(m)}$, which can be explicitly calculated:

\begin{corollary}\label{efec} The Euler characteristic of ${\rm H}_{d,n}^{(m)}$ equals $$\chi({\rm H}_{d,n}^{(m)})=\frac{n}{(m-1)d+n}{{md+n-1}\choose{d}}.$$
\end{corollary}

\proof We want to apply a formula of \cite{Sta} giving the number of plane forests.
First, we define a forest $S_*$ as above to be plane if every element $(k,w)\in S_*$ has either $m$ successors (thus $(k,wi)\in S_*$ for all $i=1,\ldots,m)$ or no successors (thus, $(k,wi)\not\in S_*$ for all $i=1,\ldots,m$). \\[1ex]
We construct a bijection $${\cal F}_{d,n}^{(m)}\leftrightarrow\{\mbox{ plane forests in }{\cal F}_{md+n,n}^{(m)}\}$$
as follows: to a forest $S_*\in{\cal F}_{d,n}^{(m)}$, associate the plane forest $T_*=S_*\cup C(S_*)$. Conversely, given a plane forest $T_*\in{\cal F}_{md+n,n}^{(m)}$, associate to it the subtree consisting of all nodes with $m$ successors. It is easy to see that this yields the claimed bijection.\\[1ex]
Now a plane forest in ${\cal F}_{md+n,n}^{(m)}$ is precisely a ``plane forest of type $$((m-1)d+n,0,\ldots,0,d)"$$ in the sense of \cite[Section 5.3]{Sta}, where the entry $d$ is placed in position $m$. By \cite[Theorem 5.3.10]{Sta}, the number of such forests is precisely
$$\frac{n}{md+n}{{md+n}\choose{d}},$$
which coincides with the claimed formula by an easy calculation. \hb

\section{Properties of the generating functions}\label{properties}

In this section, we will consider the formal power series
$$\zeta^{(m)}_n(q,t)=\sum_{d=0}^\infty\sum_k\dim A_k({\rm H}_{d,n}^{(m)})q^kt^d\in{\bf Q}[q][[t]].$$
First we give a heuristic motivation for the terminology "zeta function":\\[1ex]
Let $C$ be a smooth affine {\it curve} over the finite field ${\bf F}_q$ with $q$ elements.
Then the Hasse-Weil zeta function of $C$ is defined by
$$\zeta_C^{HW}(q,t)=\exp \sum_{d=1}^\infty |C({\bf F}_{q^d})|\frac{t^d}{d}.$$
Define also the motivic zeta function (see \cite{Kapranov}) of $C$ by
$$\zeta_C^{mot}(q,t)=\sum_{d=0}^\infty |C^{(d)}({\bf F}_q)|t^d,$$
where $C^{(d)}=C^d/S_d$ denotes the $d$-th symmetric product of $C$.
Then $$\zeta^{HW}_C=\zeta^{mot}_C$$ by \cite{DeLo}. On the other hand, we have an isomorphism
$C^{(d)}\simeq{\rm Hilb}^d(C)$ to the Hilbert scheme of $d$ points in $C$,
so that the motivic zeta function of $C$ can also be viewed as the generating function
for the number of ideals of finite codimension in the coordinate ring ${\bf F}_q[C]$.\\[1ex]
On the other hand, we can argue that none of the constructions in section \ref{celldec} depended on the ground field, which yields
$$\zeta^{(m)}_n(q,t)=\sum_{d=0}^\infty|{\rm H}_{d,n}^{(m)}({\bf F}_q)|t^d,$$
and thus $\zeta^{(m)}_n(q,t)$ can be viewed as the generating function for the number of ideals of finite codimension in $A$.\\[1ex]
By M. Kontsevich's philosophy of non-commutative geometry \cite{Kon}, the free algebra $A$ has some {\it curve}-like behaviour, since ${\rm gl}\dim A=1$.\\[1ex]
Thus, it is tempting to view $\zeta^{(m)}_1(q,t)$ as the zeta function of ``the non-com\-mu\-ta\-tive variety with coordinate ring $A$". This makes it desirable to study properties of the functions $\zeta_n^{(m)}(q,t)$ in general.\\[2ex]
Using the cell decomposition \ref{cd}, this immediately reduces to a combinatorial problem: by Corollary \ref{it}, we have
$$\zeta^{(m)}_n(q,t)=\sum_{S_*}q^{d(S_*)}t^{|S_*|},$$
where the sum runs over the set
$${\cal F}_n^{(m)}=\bigcup_d{\cal F}_{d,n}^{(m)}$$
of all $m$-ary forests with $n$ roots.\\[1ex]
It turns out that a slightly modified form of this zeta-function is more convenient. 
Therefore, we define $$d'(S_*)=d(S_*)-(m-1)\frac{|S_*|(|S_*|+1)}{2}-|S_*|.$$
This has the effect of ``normalizing" the Poincar\'e polynomial in intersection theory, so that it has a non-zero value at $q=0$. We introduce the modified zeta-function
$$\overline{\zeta}^{(m)}_n(q,t)=\sum_{S\in{\cal F}_n^{(m)}}q^{d'(S_*)}t^{|S_*|}.$$

{\bf Remark:} No nice structural properties, like for example functional equations, have been found for the unmodified generating function $\zeta_{n}^{(m)}$. It would be interesting to have a conceptual explanation for this.\\[2ex]
To study the modified zeta-function, we first restrict to the case $n=1$, thus we consider $m$-ary trees with $d$ nodes. On such trees, we have the operation of grafting. Translated to our description of trees via subsets $S$ of $\Omega$ stable under taking left subwords, we have the following:

\begin{definition} Given trees $S_1,\ldots,S_m\in\Omega$, define the grafting of $S_1,\ldots, S_m$ as the tree $$g(S_1,\ldots,S_m)=\{iw\, :\, w\in S_i,\, i=1,\ldots,m\}\cup\{()\}.$$
\end{definition}

Note that the cardinality $|g(S_1,\ldots,S_m)|$ of the grafting equals $\sum_{i=1}^m|S_i|+1$.\\[1ex]
The grafting operation clearly gives a bijection
$$g:\bigcup_{d_1+\ldots+d_m=d-1}{\cal F}_{d_1,1}^{(m)}\times\ldots\times{\cal F}_{d_m,1}^{(m)}\stackrel{\sim}{\rightarrow}{\cal F}_{d,1}^{(m)}$$
for all $d\geq 1$. We aim at a description of $d(g(S_1,\ldots,S_m))$.

\begin{lemma} For any $m$-tuple of trees $S_1,\ldots,S_m$, there is a bijection
$$D(g(S_1,\ldots,S_m))\simeq \bigcup_{i=1}^mD(S_i)\cup\bigcup_{1\leq i<j\leq m}S_i\times C(S_j)\cup\bigcup_{i=1}^mC(S_i).$$
\end{lemma}

\proof Suppose we are given $(w_1,w_2)\in D(g(S_1,\ldots,S_m))$. If $w_1=()$, then $w_2$ is of the form $w_2=iw_2'$ for a word $w_2'\in S_i$. This gives the third component in the claimed bijection.\\[1ex]
 If $w_1\not=()$, then clearly $w_2\not=()$, and thus we can write $$w_1=iw_1'\mbox{ and }w_2=jw_2'$$ for words $w_1'\in S_i$ and $w_2'\in C(S_j)$. Note that $i\leq j$. If $i<j$, then there is no further condition on $w_1',w_2'$, giving the second component in the claimed bijection.\\[1ex]
If $i=j$, then clearly $(w_1',w_2')$ has to belong to $D(S_i)$, giving the first component. \hb

\begin{corollary}\label{dg} For any $m$-tuple of trees $S_1,\ldots,S_m$, we have
$$d(g(S_1,\ldots,S_m))=\sum_{i=1}^md(S_i)+(m-1)\cdot\sum_{i<j}|S_i|\cdot|S_j|+
\sum_{i=1}^m(2m-i-1)|S_i|+m.$$
\end{corollary}

\proof This follows immediately from the previous lemma, together with the formula for $c(S_i)$. \hb

\begin{corollary} For any $m$-tuple of trees $S_1,\ldots,S_m$, we have
$$d'(g(S_1,\ldots,S_m))=\sum_{i=1}^md'(S_i)+\sum_{i=1}^m(m-i)|S_i|.$$
\end{corollary}

\proof This is a direct calculation using the previous corollary, and the formula $|g(S_1,\ldots,S_m)|=\sum_{i=1}^m|S_i|+1$. \hb

\begin{theorem}\label{fe1} The modified zeta function $\overline{\zeta}_1^{(m)}$ fulfills the functional equation
$$\overline{\zeta}^{(m)}_1(q,t)=1+t\cdot\prod_{i=0}^{m-1}\overline{\zeta}^{(m)}_1(q,q^it).$$
\end{theorem}

\proof We use the fact that any non-empty tree can be written as a grafting, together with the previous corollary, to calculate:

\begin{eqnarray*}
\overline{\zeta}^{(m)}_1(q,t)&=&\sum_{S\in{{\cal F}_1^{(m)}}}q^{d'(S)}t^{|S|}\\
&=&1+\sum_{S_1,\ldots,S_m\in{{\cal F}_1^{(m)}}}q^{d'(g(S_1,\ldots,S_m))}t^{|g(S_1,\ldots,S_m)|}\\
&=&1+t\sum_{S_1,\ldots,S_m}q^{\sum_{i=1}^md'(S_i)}q^{\sum_{i=1}^m(m-i)|S_i|}
t^{\sum_{i=1}^m|S_i|}\\
&=&1+t\sum_{S_1,\ldots,S_m}q^{\sum_{i=1}^md'(S_i)}(q^{m-1}t)^{|S_1|}(q^{m-2}t)^{|S_2|}\cdots(q^{m-m}t)^{|S_m|}\\
&=&1+t\prod_{i=1}^m\overline{\zeta}^{(m)}_1(q,q^{m-i}t).
\end{eqnarray*}\hb

It is now an easy task to generalize from the case $n=1$ to arbitrary $n$. 

\begin{lemma} For a forest $S_*\in{\cal F}_{d,n}^{(m)}$, we have
$$d'(S_*)=\sum_{i=1}^nd'(S_i)+\sum_{i=1}^n(n-i)|S_i|.$$
\end{lemma}

\proof This is proved similarly to Corollary \ref{dg}. \hb

Finally, similarly to Theorem \ref{fe1}, we prove:

\begin{proposition}\label{fe15} The modified zeta-function $\overline{\zeta}^{(m)}_{n}$ is given by
$$\overline{\zeta}_n^{(m)}(q,t)=\prod_{i=0}^{n-1}\overline{\zeta}_1^{(m)}(q,q^it).$$
\end{proposition}

\proof (of Theorem \ref{main2}): Using Corollary \ref{betti}, the function $\overline{\zeta}_n^{(m)}(q,t)$ can be identified with the function defined in the statement of Theorem \ref{main2}. The functional equations are precisely Theorem \ref{fe1} and Proposition \ref{fe15}. Uniqueness is immediate.\hb

\begin{corollary}\label{fe2} The specialization $\zeta^{(m)}_n(1,t)\in{\bf Q}[[t]]$ at $q=1$, which is the generating function for the Euler characteristics of the ${\rm H}_{d,n}^{(m)}$, is uniquely determined by the functional equations
$$\zeta^{(m)}_n(1,t)=\zeta^{(m)}_1(1,t)^n,\;\;\; \zeta^{(m)}_1(1,t)=1+t\cdot\zeta^{(m)}_1(1,t)^m.$$
\end{corollary}

\examples
\begin{itemize}
\item In the case $m=1$, we have
$$\overline{\zeta}_n^{(1)}(q,t)=\frac{1}{(1-t)\cdot\ldots\cdot(1-q^{n-1}t)},$$
as already observed in \cite{LR}.
\item In the case $m=2$, we have
$$\zeta_1^{(2)}(1,t)=\frac{1-\sqrt{1-4t}}{2t},$$
which is the well-known generating function for the Catalan numbers.
\item The modified zeta function $\overline{\zeta}_1^{(2)}(q,t)$ itself at least admits a continued fraction expansion
$$\overline{\zeta}_1^{(2)}(q,t)=\frac{1}{1-\frac{t}{1-\frac{qt}{1-\frac{q^2t}{1-\ldots}}}}.$$
\item It is easy to see that the function $\zeta_1^{(m)}(1,t)$ has a singularity at $t=\frac{(m-1)^{(m-1)}}{m^m}$, with value $\frac{m}{m-1}$. This, together with Corollary \ref{fe2}, yields the explicit formula
$$\sum_{d=0}^\infty\chi({\rm H}_{d,n}^{(m)})\cdot\left(\frac{(m-1)^{(m-1)}}{m^m}\right)^d=\left(\frac{m}{m-1}\right)^n.$$
\end{itemize}

As a conclusion to the discussion of the beginning of that section, we see that the (modified) zeta-functions $\overline{\zeta}_1^{(m)}(q,t)$ of the free associative algebras satisfy some algebraic functional equation, although they are not rational for $m\geq 2$ (which can be proven by first assuming rationality: $\zeta^{(m)}_1(1,t)=R(t)/S(t)$ for some polynomials $R(t),S(t)\in{\bf Q}[t]$, applying the functional equation \ref{fe2}, and computing the degree of $R(t),S(t)$ to get a contradiction). \\[2
ex]
We finish by giving (without proof) three other properties of $\zeta^{(m)}_n(q,t)$:
\begin{itemize}
\item The relation between trees and lattice paths (see \cite[Proposition 6.2.1]{Sta}) can be extended to provide a bijection between ${\cal F}_{d,n}^{(m)}$ and the set of lattice paths from $(0,0)$ to $(d,(m-1)d+n-1)$, using steps $(1,0)$ and $(0,1)$, and never rising above the line defined by $y=(m-1)x+n-1$. The statistics $d'$ on forests then translates into the area statistics (see \cite{Tac}) for such lattice paths. Using this, one can prove that
$$\zeta^{(m)}_n(q,t)=\sum_{d=0}^\infty\sum_\lambda q^{(m-1)\frac{d(d-1)}{2}+(n-1)d-|\lambda|}t^d,$$
where the inner sums runs over all sequences $0\leq\lambda_1\leq\ldots\leq\lambda_d$ of integers such that
$$\lambda_i\leq(m-1)i+n-1\mbox{ for all }i=1,\ldots,d.$$
\item Consider the $q$-hypergeometric function
$$\gamma^{(m)}(q,t)=\sum_{d=0}^\infty(-1)^d\frac{q^{(m-1)\frac{d(d-1)}{2}}}{(1-q)\cdots(1-q^d)}t^d.$$
Then
$$\overline{\zeta}^{(m)}_n(q,t)=\frac{\gamma^{(m)}(q,q^nt)}{\gamma^{(m)}(q,t)}.$$
This is an easy generalization of \cite[Solution to Exercise 6.34.]{Sta}; see also \cite{Duchon1}.

\item The modified zeta-function $\overline{\zeta}^{(2)}_1(q,t)$ appears in the context of {\it commutative} Hilbert schemes of points in the plane ${\bf A}^2$, see \cite{Haiman}. It would be very interesting to find an explanation.
\end{itemize}

\section{Asymptotics}\label{asymp}

The explicit formula \ref{efec} allows us to consider the asymptotic behaviour of the Euler characteristic $\chi({\rm H}_{d,n}^{(m)})$ for large $d$. Using the Stirling approximation $$n!\sim \sqrt{2\pi}\cdot n^{n+\frac{1}{2}}\cdot e^{-n},$$ one easily derives:

\begin{proposition} The asymptotic behaviour of the Euler characteristic of ${\rm H}_{d,n}^{(m)}$ for large $d$ is given by:
$$\chi({\rm H}_{d,n}^{(m)})\sim \left(\frac{n}{\sqrt{2\pi}}\cdot\left(\frac{m}{m-1}\right)^{n+\frac{1}{2}}\right)\cdot d^{-\frac{3}{2}}\cdot \left(\frac{m^m}{(m-1)^{(m-1)}}\right)^d.$$
\end{proposition}

This shows a fundamental difference to the asymptotic behaviour of the Euler characteristic of the Hilbert schemes ${\rm Hilb}^d({\bf A}^2)$ of length $d$ subschemes of the affine plane, which equals the number $p(d)$ of partitions of $d$ by \cite{Goettsche}. By the classical result of Hardy-Ramanujan, the asymptotic behaviour is thus $$\chi({\rm Hilb}^d({\bf A}^2))\sim \frac{1}{4\sqrt{3}}\cdot d^{-1}\cdot\exp(\pi\sqrt{2/3\cdot d}).$$
Note that the behaviour is thus sub-exponential in this case, whereas it is exponential in the case of the non-commutative Hilbert schemes.\\[2ex]
Finally, it is possible to identify the Airy distribution as a limit distribution of the Betti numbers of the varieties ${\rm H}_{d,1}^{(m)}$.\\[1ex]
The Airy distribution is defined as the law of a positive random variable $X$ which is uniquely determined by its $k$-th moments
$${\bf E}(X^k)=\frac{2\sqrt{\pi}}{\Gamma(\frac{3k-1}{2})}\cdot\Omega_k,$$
where the $\Omega_k$ are defined recursively by $\Omega_{0}=-1$ and
$$2\Omega_k=(3k-4)k\Omega_{k-1}+\sum_{i=1}^{k-1}{k\choose i}\Omega_i\Omega_{k-i}.$$
See \cite{FL} for a detailed discussion of this distribution.\\[1ex]
Using the results of section \ref{properties}, we can easily apply a theorem due to P.~Duchon \cite[Theorem 2]{Duchon} to get:

\begin{theorem}\label{airy} For each $d\in{\bf N}$, define a discrete random variable $X_d$ by
$${\bf P}(X_d=k)=\frac{1}{\chi({\rm H}_{d,1}^{(m)})}\cdot\dim H^{(m-1)d(d-1)-2k}({\rm H}_{d,1}^{(m)}).$$
Then the sequence of discrete random variables
$$\sqrt{\frac{8}{m(m-1)}}\cdot d^{-\frac{3}{2}}\cdot X_d$$
has the Airy distribution as a limit law.
\end{theorem}

\end{document}